\documentclass[11pt,english]{article}
\usepackage[english]{babel}
\usepackage{amsfonts}
\usepackage{epsfig}
\usepackage{graphicx}
\usepackage{fancyhdr}
\usepackage{amssymb}
\usepackage{amsmath}
  \usepackage{verbatim} 
  \usepackage{appendix}

\newtheorem{theorem}{Theorem}

\newtheorem{proposition}[theorem]{Proposition}
\newtheorem{prop}{Proposition}[section]
\newtheorem{lemma}[prop]{Lemma}
\newtheorem{rem}{Remark}[section]
\newtheorem{definition}[prop]{Definition}
\newtheorem{cl}[prop]{Claim}

\def\Lg {{\cal L}} 

\def\and {{\rm \; and \;}}

\def\exp {{\rm exp}}

\def\Lg {{\cal L}} 

\def\and {{\rm \; and \;}}

\def\exp {{\rm exp}}

\def\T0{T_{0,1}}

\def\Kt{K_{01}}
\def\det{\delta_{1}}
\def\At{A_1}
\def\st{s_{0,0}}

\def\Kc{K_{02}}
\def\sc {s_{0,1}}

\def\der{\delta_2}
\def\Ar{A_2}
\def\sr{s_{0,2}}

\def\str{s_{0,3}}
\def\Atr{A_3}

\def\Kd{K_{05}}
\def\sd {s_{0,4}}

\def\Kdi{K_{06}}
\def\sdi{s_{0,5}}

\def\Acn{A_4}
\def\scn{s_{0,6}}

\def\Adu{A_5}
\def\sdu{s_{0,7}}

\def\sdc{s_{0,8}}

\def\sfl{s_{0,9}}

\def\spi{s_{1,0}}

\def\Acp{A_7}

\def\tt {{\varrho}}
\def\etaa{{\varsigma}}

\newcommand {\R}{ \mathbb{R}}
\newcommand {\C}{ \mathbb{C}}
\newcommand {\N}{ \mathbb{N}}
\newcommand {\Z}{ \mathbb{Z}}
\newcommand{\dd}{d_0, d_1}
\newcommand {\pa}{\partial}

\newcommand {\beqna} {\begin{eqnarray}}
\newcommand {\eeqna} {\end{eqnarray}}
\newcommand {\beqtn} {\begin{equation}}
\newcommand {\eeqtn} {\end{equation}}

\newcommand {\pone}{${\cal{R}}_1$}
\newcommand {\ptwo}{${\cal{R}}_2$}

\begin{document}

\title{\bf Construction of a stable periodic solution to a semilinear heat equation with a prescribed profile}
\author{Fethi Mahmoudi\footnote{This author has been supported  by
Fondecyt Grant 1140311, Fondo Basal CMM and  partially by
``Millennium Nucleus Center for Analysis of PDE NC130017". This work started when F. Mahmoudi was visiting University of Paris 13 and University of Paris Dauphine. He is grateful for their kind hospitalty.},\\
{\it \small fmahmoudi@dim.uchile.cl}\\
{\it \small Centro de Modelamiento Matem\'{a}tico, Universidad de Chile, }\\
Nejla Nouaili,\\ {\it \small Nejla.Nouaili@dauphine.fr}\\ {\it \small CEREMADE, Universit\'e Paris Dauphine, Paris Sciences et Lettres.}\\ 
Hatem Zaag\footnote{This author has been supported by the ERC Advanced Grant  no. 291214, BLOWDISOL, and by ANR project no. ANR-13-BS01-0010-03, ANA\'E.},\\
{\it \small Hatem.Zaag@univ-paris13.fr}\\{\it \small Universit\'e Paris 13, Sorbonne Paris Cit\'e,}\\ {\it \small LAGA, CNRS (UMR 7539), F-93430, Villetaneuse, France.}}

\maketitle
\begin{abstract}

We construct a periodic solution to the semilinear heat equation with power nonlinearity, in one space dimension, which blows up in finite time $T$ only at one blow-up point. We also give a sharp description of its blow-up profile. The proof relies on the reduction of the problem to a finite dimensional one  and the use of index theory to conclude. Thanks to the geometrical interpretation of the finite-dimensional parameters in terms of the blow-up time and blow-up point, we derive the stability of the constructed solution with respect to initial data. 

\end{abstract}

\section{Introduction}
The behavior of partial differential equations (PDEs) may depend on the space $\Omega$ where they are considered. As a matter of fact, it is classical in the litterature to compare the case of $\Omega=\R^N$ with the case where $\Omega$ is a bounded domain, with boundary conditions. 

\medskip

However, in both cases, the space is ``flat'', which may prevent from seeing the effect of curvature, very important in many physical situations. Thus, considering the case of a new flat $\Omega$ appears to be relevant, at least for this reason, but not only. Indeed, that case appears also to be very challenging, from a mathematical point of view.

\medskip

As a matter of fact, many authors considered this case in the literature, in recent years. Let us mention the case of the nonlinear Schr\"odinger equation on $\mathbb{S}^N$ by Pausader, Tzvetkov and Wang in \cite{PTWAIHP14} and also by Burq, G\'erard and Tzvetkov in \cite{BGTMRL02}. There is also the work of M\'ehats and G\'erard \cite{GMSJMA11} who consider the Schr\"odinger-Poisson system on $\mathbb{S}^2$. We cite also the case of the Navier-Stokes equation on the sphere  $\mathbb{S}^2$ by  Cao, Rammaha and Titi \cite{CRTZAMP99}.

\medskip

In the case of the nonlinear heat equation, we mention  the work by Matano, Punzo and Tesei in \cite{MPT15}, where they considered front propagation on the hyperbolic space $\mathbb{H}^N$, $N\geq 2$. We also refer to the work by Cao, Rammaha and Titi \cite{CRTJDDE00}, where they considered the nonlinear parabolic equation on $\mathbb{S}^2$.

\medskip

In this paper, we would like to see whether the curvature may influence the blow-up behavior, by considering the following semilinear heat equation.

\beqtn
 \pa_t u=\Delta_{\mathbb{S}^N} u+|u|^{p-1}u,
\label{uv}
\eeqtn
where $u(t):x\in \mathbb{S}^N\to\R$ and $\Delta_{\mathbb{S}^N}$ denotes the Laplace-Beltrami operator on  $\mathbb{S}^N\subset \R^{N+1}$.\\

\medskip
By using the localization in charts, we felt that the interaction between the singular regions (where blow-up occurs) and the regular region highly challenging, and no easy solution seems to be available. For that reason, we restrict ourselves to the case where $N=1$ which is easier, due to the one-dimensional character and the periodicity. Of course the case $N=1$ conserves all the challenging character of the interaction question.

\medskip

The Cauchy problem for system (\ref{uv}) can be solved in $(L^\infty(\mathbb{S}^N))$, locally in time. We say that $u(t)$ blows up in finite time $T<\infty$, if $u(t)$ exists for all $t\in [0,T)$ and 
$\lim_{t\to T}\|u(t)\|_{L^\infty}=+\infty.$
In that case, $T$ is called the blow-up time of the solution. A point $x_0\in\mathbb{S}^N$ is said to be a blow-up point if there is a sequence $\{(x_j,t_j)\}$, such that $x_j\to x_0$, $t_j\to T$ and $|u(x_j,t_j)|\to \infty$ as $j\to\infty$.

\bigskip

\noindent Two questions arise in the blow-up study:
\begin{itemize}
\item the classification for arbitrary blow-up solutions,
\item the construction of examples obeying same prescribed behavior.
\end{itemize}
Let us first briefly mention some literature on blow-up in the case where $\Omega=\R^N$.  Consider the equation
\beqtn
 \pa_t u=\Delta u +|u|^{p-1}u,\;\; x\in\R^N.
\label{uscalar}
\eeqtn
The blow-up question for equation (\ref{uscalar}), has been studied intensively by many authors and no list can be exhaustive.\\
Regarding the classification question, when it comes to deriving the blow-up profile, the situation is completely understood in one space dimension (however, less is understood in higher dimensions, see Vel{\'a}zquez \cite{VCPDE92, VTAMS93,VINDIANA93} and Zaag \cite{ZIHP02, ZCMP02, ZMME02} for partial results). In one space dimension, given  a blow-up point $a$, we have the following alternative

\medskip

{\textit {
\begin{itemize}
\item either
\beqtn
\sup_{|x-a|\leq K\sqrt{(T-t)\log(T-t)}}\left|(T-t)^{\frac{1}{p-1}}u(x,t)-f\left(\frac{x-a}{\sqrt{(T-t)|\log(T-t)}|}\right)\right|\to 0,
\label{profileu}
\eeqtn
\item or for some $m\in \N$, $m\ge 2$, and $C_m>0$
 \beqtn
\sup_{|x-a|<K(T-t)^{1/2m}}\left|(T-t)^{\frac{1}{p-1}}u(x,t)-f_m\left(\frac{C_m (x-a)}{(T-t)^{1/2m}}\right)\right|\to 0,
\label{um}
\eeqtn
as $t\to T$, for any $K>0$, where 
\begin{equation}
\label{deff}
f(z)=\left(p-1+\frac{(p-1)^2}{4p}z^2\right)^{-\frac{1}{p-1}}
\mbox{ and }f_m(z) =\left(p-1+|z|^2m\right)^{-\frac{1}{p-1}}.
\end{equation}
\end{itemize}
}}
Let us mention that when $\Omega$ is a subdomain of $\R^N$ with $C^{2,\alpha}$ boundary, Giga and Kohn proved in \cite{GKCPAM89} that no blow-up occurs at the boundary and many of the above results hold.

\bigskip

Regarding the construction question, from Bricmont and Kupiainen \cite{BKN94} and Herrero and Vel{\'a}zquez \cite{HVAIHP93}, we have examples of initial data leading to each of the above-mentioned scenarios. Note that (\ref{profileu}) corresponds to the fundamental mode of the harmonic oscillator in the leading order, whereas (\ref{um}) corresponds to higher modes. Moreover, Herrero and 
Vel{\'a}zquez proved the genericity of the behavior (\ref{profileu}) in one space dimension in \cite{HVADSSP92} and \cite{HVCRAS94}, and only announced the result in the higher dimensional case (the result has never been published). Note also that the stability of such a profile with respect to initial data has been proved by Fermanian Kammerer, Merle and Zaag  in \cite{FZN00} and \cite{FMZMA00}. For more results on equation (\ref{uscalar}), see \cite{BQJMOS77}, \cite{GKCPAM85}, \cite{GKIUMJ87}, \cite{GKCPAM89}, \cite{HVAIHP93}, \cite{HVCRAS94}, \cite{MMCPAM04}, \cite{MMJFA08}, \cite{MZCPAM98}, \cite{MZMA00}, \cite{MMA07} and \cite{QSBV07}.

\bigskip

In this paper, we address the construction question of equation  \eqref{uv} when $N=1$, and we give the first example of a blow-up solution of equation \eqref{uv} in $\mathbb{S}=\mathbb{S}^1$,  with a full description of its blow-up profile, obeying behavior \eqref{profileu} (note that our method extends with no difficulty to the construction of an analogous solution obeying the behavior \eqref{um}; however, the proof should be even more technical). 

\medskip

Making the change of variables $x=e^{i\theta}\in \mathbb{S}$, we can rewrite \eqref{uv} as 
\beqtn
\pa_t u=\pa_{\theta}^{2} u+|u|^{p-1}u,
\label{eqtheta}
\eeqtn
where $\theta\in\R$ and $u(.,t)$ is $2\pi-$periodic. As we said above, the first main purpose of this work is to show that behavior (\ref{profileu}) does occur. More precisely, we prove the existence of a blow-up solution for equation (\ref{eqtheta}) and we give a description of its profile. This can be summarized in the following  result. 

\begin{theorem}[Existence of a blow-up solution for  (\ref{eqtheta}) with prescribed profile]

There exists $T>0$ such that equation (\ref{eqtheta}) has a solution $u(\theta,t)$ in $\mathbb{S}\times [0,T)$ such that:\\
\begin{enumerate}
\item[(i)] the solution $u$ blows up in finite time $T$ only in $2\pi\Z$; \\
\item[(ii)] there holds that for all $R>0$,
\beqtn
\sup_{\Lambda_R} \left | (T-t)^{\frac{1}{p-1}}u(\theta+2k\pi,t)-f \left(\frac{\theta+2k\pi}{\sqrt{(T-t)|\log(T-t)|}}\right)\right |\ \longrightarrow 0\mbox{ as }t\to T,
\label{profilev}
\eeqtn 
 where $\Lambda_R:=\big\{|\theta-2k\pi|\leq R\sqrt{(T-t)|\log(T-t)|},k\in\Z\big\}$ and where the function $f$ is defined by (\ref{deff}).\\
\item[(iii)] for all $\theta \not\in 2\pi\Z $, $u(\theta, t)\to u(\theta,T)$ as $t \to T$, with $u(\theta,T)\sim u^*(\theta-2k\pi)\mbox{ as }\theta \to 2k\pi\mbox{, for any }k\in\Z$, where
\beqtn
u^*(\theta)= \left[ \frac{(p-1)^2|\log\theta|}{8 p \theta^2}\right]^{-\frac{1}{p-1}}.
\label{defue}
\eeqtn
\end{enumerate}
\label{theorem1}
\end{theorem}
\begin{rem}

\begin{enumerate}
\item A natural extension of our result would be to consider the heat equation in $\R^N$ in the periodic
setting (or equivalently, in the torus $\R^N/\Z^N$ for example). We believe that our result also holds
in this setting. However, understanding the interaction between the inner and the outer parts (or the
solution in the ÒregularÓ and the Òblow-upÓ parts; see Section \ref{formulationproblem}) would be much more
difficult, though the difficulty is only technical. For that reason, we focus on one space dimension in
this paper to keep the paper within reasonable limits.
\item
In the case of equations \eqref{uv} and \eqref{uscalar}, there do exist solutions which
behave like (3) and others like (4) (see Bricmont-Kupiainen, \cite{BKN94}, for the construction
in both cases, when the equation is considered in $\R^N$). In our paper, we focus on the
construction of solutions obeying (3), but we do believe that our method can be adapted to construct
solutions obeying (4). However, we only construct a solution obeying (3), since this should be the
only stable solution.
\end{enumerate}
\end{rem}
The proof of Theorem \ref{theorem1} uses some ideas developed by  Bricmont and Kupiainen \cite{BKN94} and Merle and Zaag \cite{MZDuke97} and \cite{MZNL97} to construct a blow-up solution for the semilinear heat equation (\ref{uscalar}) obeying the behavior (\ref{profileu}). In \cite{EZSMJ11}, Ebde and Zaag use the same ideas to show the persistence of the profile (\ref{profileu}) under perturbations of equation (\ref{uv}) in the real case by lower order terms involving $u$ and $\nabla u$. See also \cite{NVTZ15} by Nguyen and Zaag, where the authors consider stronger perturbation than the one considered in \cite{EZSMJ11}.

 In \cite{MZ07}, Masmoudi and Zaag adapted that method to the following complex Ginzburg-Landau equation with no gradient structure 
\[\pa_t u=(1+i\beta)\Delta u+(1+i\delta)|u|^{p-1}u\mbox{, with $\beta, \,\delta\in \mathbb R$}.\]
Notice that the case $\beta=0$ and $\delta$ small was first considered by  Zaag, see \cite{ZAIHPANL98}. We also mention the recent work by Nouaili and Zaag in \cite{NZCPDE15} for a complex-valued equation with no gradient structure.
\\
More precisely, the proof relies on the understanding of the dynamics of the self-similar version of (\ref{uv}) (see system (\ref{qqtilde}) below) around the profile (\ref{profileu}). Moreover, we proceed in two steps:
\begin{itemize}
\item First, we reduce the question to a finite-dimensional problem: we show that it is enough to control a $(N+1)$-dimensional variable in order to control the solution (which is infinite dimensional) near the profile.
\item Second, we proceed by contradiction to solve the finite-dimensional problem and conclude using index theory.
\end{itemize}
Surprisingly enough, we would like to mention that this kind of methods has proved to be successful in various situations including hyperbolic and parabolic PDE, in particular with energy-critical exponents. This was the case for the construction of multi-solitons for the semilinear wave equation in one space dimension by C\^{o}te and Zaag \cite{CZCPAM13}, the wave maps by Rapha\"{e}l and Rodnianski \cite{RRIHES12}, the Schr\"{o}dinger maps by Merle, Rapha\"{e}l and Rodnianski \cite{MRRCRAS11}, the critical harmonic heat flow by Schweyer \cite{SJFA12} and the two-dimensional Keller-Segel equation by Rapha\"{e}l and Schweyer \cite{RSCPAM13}.

\medskip

The interpretation of the parameters of the finite dimensional problem in terms of the blow-up time and blow-up point allows us, as in \cite{MZDuke97} and \cite{MZNL97}, to derive the stability of the constructed solution as stated in the following result.

\begin{proposition}[Stability of the constructed solutions]Denote by $\hat{u}$ the solution constructed in Theorem \ref{theorem1} and by $\hat{T}$ its blow-up time. Then, there exists $\varepsilon_0>0$ such that for any initial data $u_0\in L^\infty(\mathbb{S})$, satisfying
\[\|u_0-\hat{u}(.,0)\|_{L^\infty}\leq \varepsilon_0,\]
the solution of equation \eqref{eqtheta}, with initial data $u_0$ blows up in finite time $T(u_0)$ at only one blow-up point $a(u_0)\in \mathbb{S}$.\\
Moreover, the function $u(\theta-a(u_0),.)$ satisfies the same estimates as $u$ with $\hat{T}$ replaced by $T(u_0)$.\\
Furthermore, it follows that
\[T(u_0)\to \hat{T},\;\;a(u_0)\to a_0\mbox{  as } u_0\to \hat{u}(0).\]
\end{proposition}
The proof of this stability result follows exactly as in \cite{MZDuke97} and \cite{MZNL97}. For that reason we skip it and refer the interested readers to those papers.

\medskip

We proceed in 3 sections to prove Theorem \ref{theorem1}. We first give in Section 2 an equivalent formulation of the problem in the scale of the well-known similarity variables. Section 3 is devoted to the proof of the similarity variables formulation (this is a central part in our argument). In the last section, we prove Theorem \ref{theorem1}. 
 \section{Formulation of the problem}
 \label{formulationproblem}
 We would like to find initial data  $u_0$ such that the solution $u$ of equation (\ref{uv}) blows up in time $T$ with
\begin{eqnarray}
\sup_{\Lambda_R} \left | (T-t)^{\frac{1}{p-1}}u(\theta+2k\pi,t)-f \left(\frac{\theta+2k\pi}{\sqrt{(T-t)|\log(T-t)|}}\right)\right |\ 
\longrightarrow  0\mbox{ as }t\to T
\end{eqnarray}
where  $\Lambda_R:=\big\{|\theta-2k\pi|\leq R\sqrt{(T-t)|\log(T-t)|},k\in\Z\big\}$
\medskip

This is the main estimate and the other results of Theorem \ref{theorem1} will appear as by-products of the proof (see Section 4 for the proof of all the estimates of Theorem \ref{theorem1}). From periodicity, we will consider $\theta$ in one period, which depends on the region we consider (see below in the definitions of the regular and the blow-up region).\\
First, we introduce the following cut-off function $\chi_0 \in C^{\infty}_{0}(\R,[0,1])$, 
\beqtn
\chi_0(\xi)=
\left\{
\begin{array}{lll}
1&\mbox{ if }& |\xi| \leq 1,\\[3mm]
0&\mbox{ if }& |\xi|\geq 2.
\end{array}
\right.
\label{chi0}
\eeqtn
In the following, we will divide our work in two parts; the blow-up region and the regular region.
\begin{itemize}
 \item In the regular region, we will study $\bar{u}$ defined by:
\beqtn
\bar{u}(\theta)=
\left\{
\begin{array}{lll}
u(\theta)\bar{\chi}(\theta)&\mbox{ if }& \theta\in [0,2\pi],\\[3mm]
0&\mbox{ if }& \theta \in \R\setminus  [0,2\pi],
\end{array}
\right .
\label{defbaru}
\eeqtn
where the function $\bar{\chi}$ is $2\pi-$periodic and defined for all 
$\xi\in [-\pi,\pi]$ by $$\bar{\chi}(\xi)=1-\chi_0\left(\frac{4\xi}{\varepsilon_0}\right)$$
with  $\varepsilon_0>0$ will be fixed small enough later. Then,  for all $\theta\in\R$, $\bar{u}(\theta)$ satisfies the following equation
\beqtn
\pa_t \bar{u}=\pa_{\theta}^{2}\bar{u}+|u|^{p-1}\bar{u}-2\bar{\chi}'\pa_\theta u -\bar{\chi}''u.
\label{eqbaru}
\eeqtn
\noindent We control $\bar{u}$ using classical parabolic estimates on $u$ as we will see in Proposition \ref{parabolicesti}, below.
\item In the blow-up region of $u(\theta,t)$, we make the following self-similar transformation of problem (\ref{eqtheta})
\end{itemize}
\beqtn
\begin{array}{c}
W(y,s)=(T-t)^{\frac{1}{p-1}}u(\theta,t),\\[3mm]
y=\frac{\theta}{\sqrt{T-t}}\mbox{, }s=-\log(T-t),
\end{array}
\label{chauto}
\eeqtn
then $W(y,s)$, for $y\in\R$, satisfy the following equation 
\beqtn
\pa_s W=\pa_{y}^{2}W-\frac 12 y \pa_{y}W-\frac{1}{p-1}W+|W|^{p-1}W.
\label{eqw}
\eeqtn
We note that for all $s\in\R$, $W$ is $2\pi e^{s/2}$  periodic.\\
 Let us define 
\beqtn
w(y,s)=
\left\{
\begin{array}{ll} 
W(y,s)\chi(y,s)&\mbox{ if  $|y|\leq \pi e^{s/2}$},\\[3mm]
0 &\mbox{ otherwise, }
\end{array}
\right.
\label{defw}
\eeqtn
with
\beqtn
\chi(y,s)=\chi_0\left(\frac{y e^{-s/2}}{\varepsilon_0}\right),
\label{Nchi}
\eeqtn
where $\chi_0$ is defined by (\ref{chi0}) and $\varepsilon_0$ will be fixed small enough later in the proof.\\
Then we multiply equation (\ref{eqw}) by $\chi(y,s)$ and we get
\[\chi\pa_s W=\chi\pa_{y}^{2}W-\frac 12 y \chi\pa_{y}W-\frac{1}{p-1}w+|W|^{p-1}w, \]
therefore

\beqtn
\begin{array}{l}
\pa_s w=\pa_{y}^{2}w-\frac 12 y \pa_{y}w-\frac{1}{p-1}w+|w|^{p-1}w+F(y,s),
\end{array}
\label{eqwchi}
\eeqtn
where
\beqtn
F(y,s)=
\left\{
\begin{array}{ll}
W\pa_s\chi - 2\pa_y \chi \pa_y W-W\pa_{y}^{2}\chi+\frac 1 2 yW \pa_y\chi
+|W|^{p-1}W\left(\chi-\chi^{p}\right) &\mbox{if } |y|\leq \pi e^{s/2},\\[3mm]
0 &\mbox{ otherwise, }
\end{array}
\right .
\label{F}
\eeqtn
\begin{rem} \rm  We note that $w$ is not periodic, and that equation (\ref{eqwchi}) is valid for all $y\in\R$.

\medskip

\noindent Now, let us Introduce 
\beqtn
w=\varphi+q,
\label{initialqqtilde}
\eeqtn
with
\beqtn
\varphi=f\left(\frac{y}{\sqrt{s}}\right)+\frac{\kappa}{4ps},
\label{deffi}
\eeqtn
where 
\beqtn
f(z)=\left(p-1+b z^2\right )^{-\frac{1}{p-1}}\mbox{, $\kappa=(p-1)^{-\frac{1}{p-1}}$ and } b=\frac{(p-1)^2}{4p}.
\label{defkappa}
\eeqtn
The problem is then reduced to constructing a function $q$ such that
\[\lim_{s\to\infty}\|q(y,s)\|_{L^\infty}=0,\]
and $q$ is a solution of the following equation for all $(y,s)\in \R\times\left [s_0(=-\log T),\infty\right)$,
\beqtn
\begin{array}{lll}
\displaystyle \pa_s q&=&\displaystyle({\cal{L}}+V)q+B(y,s)+R(y,s)+F(y,s),
\end{array}
\label{qqtilde}
\eeqtn
where
\beqtn
\displaystyle{\cal{L}}=\pa_{y}^{2} -\frac{1}{2}y\pa_y+1\mbox{, }V(y,s)=p\varphi(y,s)^{p-1}-\frac{p}{p-1},
\label{OperatorL}
\eeqtn
\beqtn
B(y,s)=|\varphi+q|^{p-1}(\varphi+q)-\varphi^p-p\varphi^{p-1}q\mbox{, },
\label{defBV}
\eeqtn
and
\beqtn
\begin{array}{lll}
\displaystyle R(y,s)&=&\pa_{y}^{2} \varphi-\frac{1}{2}y \pa_y \varphi-\frac{\varphi}{p-1}+\varphi^{p-1}-\pa_s \varphi,\\[3mm]
\displaystyle F(y,s)&=&H(y,s)+\pa_y G(y,s)\mbox{ with,}\\[3mm]
\displaystyle H(y,s)&=&W\left(\pa_s\chi+\pa_{y}^{2}\chi+\frac 12y\pa_y\chi\right)+|W|^{p-1}W(\chi-\chi^p),\\[3mm]
\displaystyle G(y,s)&=&-2\pa_y\chi W.
\label{HG}
\end{array}
\eeqtn
The control of $q$ near $0$ obeys two facts:
\begin{itemize}
\item {\bf Localization}: the fact that our profile $\varphi(y,s)$ dramatically changes its value from $1+\frac{1}{4s}$ in the region near $0$ to $\frac{1}{4s}$ in the region near infinity, according to a free boundary moving at the speed $\sqrt{s}$. This will require different treatments in the regions $|y|<2K_0\sqrt{s}$ and $2K_0\sqrt{s}<|y|<\frac \pi 2 e^{s/2}$ for some $K_0$ to be chosen.
\item {\bf Spectral information}: the fact that the operator $\Lg$ is selfadjoint, $B$ is quadratic in $q$ and  that
\[\|R(s)\|_{L^\infty}+\|V(s)\|_{L^{2}_{\rho}}\to 0 \mbox{ as } s\to \infty,
\]
from (\ref{initialqqtilde}) and (\ref{OperatorL}), which shows that the dynamics of equation (\ref{qqtilde}) near $0$ are driven by the spectral properties of $\Lg$. This will require a decomposition of the solution according to the spectrum of $\Lg$. Note that the operator $\Lg$ is self-adjoint in the Hilbert space
\[
L^{2}_{\rho}=\bigg\{ g\in L^{2}_{loc}(\R,\C)\mbox{,
 }\|g\|_{L^{2}_{\rho}}^{2}\equiv\int_{\R}|g|^2e^{-\frac{|y|^2}{4}}dy
<+\infty\bigg\}\mbox{ with }\rho(y)=\displaystyle\frac{e^{-\frac{|y|^2}{4}}}{(4\pi)^{1/2}}.
\]
\end{itemize}
The spectrum of $\Lg$ is explicitly given by 
\[{\rm spec}(\Lg)=\{1-\frac{m}{2}\mbox{, }m\in\N\}.\]
All the eigenvalues are simple, the eigenfunctions are dilations of Hermite's polynomial and  given by 
\beqtn
h_m(y)=\sum_{n=0}^{[\frac{m}{2}]}\frac{m!}{n!(m-2n)!}(-1)^n y^{m-2n}.
\label{hermite}
\eeqtn 
Note that $\Lg$ has two positive (or expanding) directions ($\lambda=1$ and $\lambda=\frac 1 2$), and a zero direction ($\lambda=0$).
Complying with the localization and spectral information facts, we will decompose $q$ accordingly as stated above:
\begin{itemize}
\item First, let us introduce 
\beqtn
\chi_1(y,s)=\chi_{0}\left(\frac{|y|}{K_0\sqrt{s}}\right),
\label{chi1}
\eeqtn
where $\chi_0$ is defined in (\ref{chi0}), $K_0\geq 1$ will be chosen large enough so that various technical estimates hold. Then, we write $q=q_b+q_e$, where the inner part and the outer part are given by 
\beqtn
q_b=q\chi_1\mbox{, } q_e=q(1-\chi_1)  .
\label{qbqe}
\eeqtn
Let us remark that 
\[{\rm supp}(q_b(s))\subset B(0,2K_0\sqrt{s})\mbox{ and }\quad{\rm supp}(q_e(s))\subset \R\setminus B(0,K_0\sqrt{s}).\]

\item Second, we study $q_b$  using the structure of $\Lg$, isolating the nonnegative directions. More precisely we decompose $q_b$  as follows
 \beqtn
 \begin{array}{lll}
 q_b(y,s)&=&\sum\limits_{0}^{2}q_m(s) h_m(y)+q_-(y,s),
  \label{decompositionq}
 \end{array}
 \eeqtn
 where $q_m$ is the projection of $q_b$  on $h_m$, $q_-(y,s)=P_-(q_b)$ and $P_-$ is the projection on $\{h_i,\; i\geq 3\}$ the negative subspace of the operator $\Lg$. 
\end{itemize}
In summary, we can decompose $q$ in 5 components as follows:
 \beqtn
 \begin{array}{lll}
 q(y,s)&=& \sum\limits_{m=0}^{2}q_m(s)\, h_m(y)+ q_-(y,s)+q_e(y,s).
  \end{array}
 \label{decompq}
 \eeqtn
 Here and throughout the paper, we call  $q_-(y,s)$ the negative part of $q$ and $q_2$, the null mode of $q$. 
\end{rem}
 \label{Section2}
 \section{The construction method in selfsimilar variables}
 This section is devoted to the proof of the existence of a solution $u$ of equation \eqref{uv} satisfying $\|q(s)\|_{L^\infty}\to 0$. This is a central argument in our proof. Though we refer to the earlier work by Merle and Zaag \cite{MZDuke97} for purely technical details, we insist on the fact that we can completely split from that paper as long as ideas and arguments are considered. We hope that the explanation of the strategy we give in this section will be more reader friendly.

\medskip

\noindent We proceed in 3 subsections:
\begin{itemize}
\item In the first subsection, we give all the arguments of the proof without the details, which are left for the following subsection (readers not interested in technical details may stop here).
\item In the second subsection, we give various estimates concerning initial data.
\item In the third subsection, we give the dynamics of system (\ref{qqtilde}) near the zero solution, in accordance with the decomposition (\ref{decompq}) and taking into account the interaction between the singular region and the regular region.
\end{itemize}
\subsection{The proof without technical details}

Given $T>0$, we consider initial data for equation (\ref{uv}) $2\pi-$periodic defined for all $\theta\in[-\pi,\pi]$ by:
\beqtn
u_0(\theta,d_0,d_1)=T^{-\frac{1}{p-1}}\left\{\varphi(y,s_0)\chi(8y,s_0)+\frac{A}{s_{0}^{2}}(d_0+d_1 y)\chi_1(2y,s_0) \right\},
\label{initialq}
\eeqtn
where $s_0=-\log T$, $y=\frac{\theta}{\sqrt{T}}$,  $\chi$ is defined in (\ref{Nchi}) and $\chi_1$ is defined in (\ref{chi1}).

\medskip
 
\noindent Notice that $u_0$ depends also on $K_0$, $\varepsilon_0$, $A$ and $T$, but we omit that dependence in \eqref{initialq} for simplicity.\\
Notice also that the transition at $-\pi+2k\pi$ in $u_0$ is smooth, since $u_0\equiv 0$ is some open interval around that number.\\
Thanks to Section 2, in order to control $u(s)$ near $\varphi$, it is enough to control it in some shrinking set defined as follows:
\begin{definition}[Definition of a shrinking set for the components of $q$]\label{defVA}For all $K_0>0$, $\varepsilon_0>0$, $A>0$, $0<\eta_0\leq 1$ and $T>0$, we define for all $t\in[0,T)$ the set $S^*(K_0,\varepsilon_0,A,\eta_0,T,t)$ as being the set of all functions $u\in L^\infty(\R)$ satisfying:
\begin{enumerate}
\item[\rm (i)]\textbf{Estimates in \pone:} $q(s)\in V_{K_0,A}(s)$ where $s=-\log (T-t)$, $q(s)$ is defined in (\ref{chauto}), (\ref{defw}), (\ref{initialqqtilde}) and (\ref{deffi}) and $ V_{K_0,A}(s)$ is the set of all functions $r\in L^\infty(\R)$ such that
\beqtn
\left\{
\begin{array}{ll}
|r_m(s)|\leq A s^{-2}(m=0,1),&|r_2(s)|\leq A^2 s^{-2}\log s,\\[3mm]
|r_-(y,s)|\leq A s^{-2}(1+|y|^3),&|r_e(y,s)|\leq A s^{-1/2},
\label{estR1}
\end{array}
\right.
\eeqtn
where
\beqtn
\left\{
\begin{array}{ll}
r_e(y,s)=(1-\chi_1(y,s))r(y,s),&r_-(s)=P_-(\chi_1(s)r),\\[3mm]
\mbox{for $m\in\N$, }r_m(y,s)={\displaystyle \int d\rho k_m(y) \chi_1(y,s)r(y),}&
\end{array}
\right.
\eeqtn
$\chi_1$ is defined in (\ref{chi1}), $P_-$ is the $L^{2}_{\rho}$ projector on ${\rm Vect}\{h_m|m\geq 3\}$.

\item[\rm (ii)]\textbf{Estimates in \ptwo:} For all $\frac{\varepsilon_0}{2}\leq |\theta|\leq \pi$, $|u(\theta,t)|\leq \eta_0$.
\end{enumerate}
\end{definition}
\begin{rem}For simplicity, we may write $S^*(t)$ instead of $S^*(K_0,\varepsilon_0,A,\eta_0,T,t)$. Note also that our arguments work with $\eta_0=1$.\end{rem}
Our aim becomes then to prove  the following result.
\begin{prop} [Existence of a solution of (\ref{qqtilde}) trapped in $S^*(t)$] \label{Newprop}

\

There exists $\Kt>0$ such that for each $K_0\geq \Kt$, there exists $\det(K_0)$, such that for any $\varepsilon_0\leq\det(K_0)$, there exists $\At(K_0,\varepsilon_0)$, such that for any $A\geq \At$ and $0<\eta_0\leq 1$, there exists $\st(K_0,\varepsilon_0,A,\eta_0)$ such that for all $T\leq e^{-\st}$, there exists $(d_0,d_1)\in \R^2$, such that:\\
if $u$ is a solution of (\ref{eqtheta}) with initial data given by (\ref{initialq}), then 
\[\forall t\in [0,T),\; u(t)\in S^*(K_0,\varepsilon_0,A,\eta_0,T,t).\] 
\end{prop}
To prove this proposition we need some intermediate lemmas.

\medskip

\noindent In the following lemma, we find a set $D_{K_0,\varepsilon_0,A,T}=D_T$ such that $u(0)\in  S^*(0)$, whenever $(\dd)\in D_{T}$. More precisely, we claim the following: 

\begin{lemma}[Choice of parameters $\dd$ to have initial data in $S^*(0)$] \label{initialisationN}  There exists $\Kc>0$ such that for each $K_0\geq \Kc$ there exists $\varepsilon_0>0$, $A\geq 1$, there exists $\sc (K_0,\varepsilon_0,A)\geq 0$ such that for all $s\geq \sc$:\\
If initial data for equation (\ref{uv}) are given by (\ref{initialq}):
 then, there exists a rectangle 
\beqtn
D_{K_0,\varepsilon_0,A,T}=D_{T} \subset [-2,2]^2,
\eeqtn
 such that, for all $(\dd)\in D_{T}$, we have 
\[u(K_0,T,A,\dd)\in S^*(0).\]
\end{lemma}

\textit{Proof}: The proof is purely technical and follows as the analogous step in \cite{MZDuke97}, for that reason we refer the reader to Lemma 3.5 page 156 and Lemma 3.9 page 160 in \cite{MZDuke97}. $\blacksquare$

\medskip

\noindent Let us consider $(\dd)\in D_{T}$ and $s_0=-\log T\geq s_{0,1}$ defined in Lemma  \ref{initialisationN}. Since $u_0$ is $2\pi-$periodic, from the local Cauchy theory, we define a maximal $2\pi-$periodic solution $u$ to equation (\ref{uv}) with initial data (\ref{initialq}), and a maximal time $t_*(\dd)\in [0,T)$ such that, 
\beqtn
\mbox{for all }t\in[0,t_*),\; u(t)\in S^*(t), 
\label{us}
\eeqtn
\begin{itemize}
\item either $t_*=T$,
\item or $t_*< T$ and from continuity, 
\beqtn
u(t_*)\in\pa S^*(t_*), 
\label{336}
\eeqtn
in the sense that when $t=t_*$, one '$\leq$' symbol in the definition of $S^*(t_*)$ is replaced by the symbol '$=$'. 
\end{itemize}
Our aim is to show that for all $A$ and $T$ small enough, one can find a parameter $(\dd)$ in $D_{T}$ such that
\beqtn
t_*(\dd)=T.
\label{butN}
\eeqtn
We argue by contradiction, and assume that for all $(\dd) \in D_{T}$, $t_*(\dd)<T$. 
As we have just stated, one of the symbols `$\leq$' in the definition of $S^*(t)$ should be replaced by `$=$' symbols when $t=t_*$.\\ 
In fact, no `$=$' sign occurs for $q_2$, $q_{-}$, $q_e$ and the estimate in \ptwo, as one sees in the following lemma.
\begin{lemma}[Reduction to a finite dimensional problem]\label{transversality}
For any $K_0>0$, there exists $\der(K_0)>0$, such that for any $\varepsilon_0\leq\der(K_0)$, there exists $\Ar(K_0,\varepsilon_0)$, such that for $A\geq \Ar$ and $0<\eta_0\leq 1$, there exists $\sr(K_0,\varepsilon_0,A,\eta_0)$ such that for any $s\geq \sr$, we have 
\[(q_0(s_*),q_1(s_*))\in \pa\left( \left[-\frac{A}{s_{*}^{2}},\frac{A}{s_{*}^{2}}\right]^2\right)\mbox{, where }
s_*=-\log(T-t_*).\]
\end{lemma}
\begin{rem} The choice of parameters cited below is particularly intricate. We explain that at conclusion of Part 1 and Part 2 in pages \pageref{pagepart1} and \pageref{pagepart2} below.\end{rem}
\textit{Proof}: This is a direct consequence of the dynamics of equation (\ref{qqtilde}), as we will show in Subsection \ref{reduction} below.\\
Just to give a flavor of the argument, we invite the reader to look at Proposition \ref{prop36} below, where we project equation (\ref{qqtilde}) on the different components of $q$ introduced in (\ref{decompq}). There, one can see that the components $q_2$, $q_-$ and $q_e$ correspond to decreasing directions of the flow and since they are ``small'' at $s=s_0=-\log T$
 (see Lemma \ref{prop37} below), they remain small for $s\in [s_0,s_*]$, and cannot touch the boundary of the intervals imposed by the definition of $S^*(t)$ in \eqref{estR1}. Thus, only $q_0$ or $q_1$ may touch the boundary of the intervals in \eqref{estR1} at $s=s_*$.\\
For more details on the arguments, see Subsection \ref{reduction} below. This ends the proof of  Lemma \ref{transversality}.$\blacksquare$

\medskip

\noindent From Lemma \ref{transversality}, we may define the rescaled flow $\Phi$ at $s=s_*$ for the two expanding directions, namely $q_0$ and $q_1$, as follows:
\beqtn
\begin{array}{lccl}
\Phi : & D_{T}&\longrightarrow& \pa([-1,1]^2)\\
             &(d_0,d_1)&\longmapsto& \displaystyle \left(\frac{s_{*}^{2}q_0}{A},\frac{s_{*}^{2}q_1}{A}\right)_{d_0,d_1}(s_*) .    
\end{array}
\label{defphi}
\eeqtn
In particular, 
\beqtn
\mbox{either }\,\omega q_0(s_*)=\frac{A}{s_{*}^{2}}, \quad\hbox{or } \quad \omega q_1(s_*)=\frac{A}{s_{*}^{2}}
\label{qmqtildem}
\eeqtn
and  $\omega\in \{-1,1\}$, both depending on
$(\dd)$. 
In the following lemma, we show that $q_m$ actually crosses its boundary at $s=s_*$, resulting in the continuity of $s_*$ and $\Phi$. More precisely, we have the following result.
\begin{lemma}[Transverse crossing] \label{transcross}
For each $K_0>0$, $\varepsilon_0>0$, there exists $\Atr(K_0,\varepsilon_0)>0$ such that for all $A\geq \Atr$, there exists $\str(K_0,\varepsilon_0,A)$ such that if $s_0 \geq \str$, $0<\eta_0\leq1$ and (\ref{qmqtildem}) holds, then 
 \beqtn
\omega\frac{d q_{m}}{ds}(s_*)>0
\label{qmqtildemin}
\eeqtn
\end{lemma}
Clearly, from the transverse crossing, we see that
 \[
(\dd)\mapsto s_*(\dd) 
\mbox{ is continuous,}
\]
 hence by definition (\ref{defphi}), $\Phi$ is continuous.
In order to find a contradiction and conclude, we crucially use the particular form we choose for initial data in (\ref{initialq}). More precisely, we have the following:
\begin{lemma}[Degree $1$ on the boundary]\label{lem3}
There exists $\Kd>0$ such that for each $K_0\geq \Kd$, $\varepsilon_0>0$ and $A\geq 1$, there exists $\sd(K_0,\varepsilon_0,A)$ such that if $s_0\geq \sd$, then the mapping $(d_0,d_1)\to (q_0(s_0),q_1(s_0))$ maps $\pa D_{T}$ into $\pa\left ( [-\frac{A}{s_{0}^{2}},\frac{A}{s_{0}^{2}}]^2\right)$, and has degree one on the boundary.
\end{lemma}
Indeed, from this lemma and the transverse crossing property of Lemma \ref{transcross}, we see that if $(\dd)\in\pa D_{T}$, then
 $s_*(\dd)=s_0$,  

$$\Phi(s_*(s,\dd),\dd)=\left(\frac{s^{2}_{*}q_0}{A},\frac{s^{2}_{*}q_1}{A}\right)(s_0)$$
 and $\Phi$ defined in (\ref{defphi}) is a continuous function from the rectangle $D_{T}\subset \R^2$ to $\pa [-1,1]^2$, whose restriction to $\pa D_{T}$ is of degree 1. This is a contradiction.\\ 
Thus, there exists $d_0,d_1\in D_T$ such that 
\beqtn
t_*(d_0,d_1)=T\mbox{ and }\forall t\in [0,T),\;\; u(t)\in S^*(t),
\label{conpro}
\eeqtn
and Proposition \ref{Newprop} follows at once.

\medskip

\noindent Thus, Proposition \ref{Newprop} is proved, and from (i)  of Definition \ref{defVA}, we have constructed a solution $q$ to system (\ref{qqtilde}), such that
\[\|q(s)\|_{L^\infty}\to 0\mbox{ as }s\to\infty.\]
The next subsections will be devoted to the proofs of the technical  Lemmas \ref{initialisationN}, \ref{transversality}, \ref{transcross} and \ref{lem3}, referring to earlier works when the proof is the same.
 \subsection{Preparation of initial data}
In this subsection, we study initial data given by (\ref{initialq}). More precisely, we state a lemma which directly implies Lemmas \ref{initialisationN} and \ref{lem3}. It also shows the (relative) smallness of the components $q_2$, $q_-$ and $q_e$, an information which will be useful for the next subsection, dedicated to the dynamics of equation (\ref{qqtilde}), crucial for the proofs of the reduction to a finite dimensional problem (Lemma \ref{transversality}) and the transverse crossing property (Lemma \ref{transcross}). More precisely, we claim the following:

\begin{lemma}[Decomposition of initial data in different components]\label{prop37}
There exists $\Kdi>0$ such that for each $K_0\geq \Kdi$, $\varepsilon_0>0$, $A\geq 1$, there exists 
$\sdi(K_0,\varepsilon_0A)\geq e$ 
such that for all $s_0\geq \sdi$
\begin{enumerate}
\item[\rm(i)] there exists a rectangle
\beqtn
D_{K_0,\varepsilon_0,A,T}=D_{T} \subset [-2,2]^2,
\eeqtn
such that the mapping $(d_0,d_1)\to (q_0(s_0),q_1(s_0))$ is linear and one to one from $D_{T}$ onto $[-\frac{A}{s_{0}^{2}},\frac{A}{s_{0}^{2}}]^2$ and maps the boundary $\pa D_{T}$ into the boundary $\pa\left ( [-\frac{A}{s_{0}^{2}},\frac{A}{s_{0}^{2}}]^2\right)$. Moreover, it is of degree one on the boundary.
\item[\rm(ii)] For all $(d_0,d_1)\in D_{T}$, we have:
\beqtn
 \begin{array}{l}
 |q_2(s_0)|\leq C A e^{-s_0},\;\; |q_-(y,s_0)|\leq \frac{c}{s_{0}^{2}}(1+|y|^3)\mbox{ and }q_e(y,s_0)= 0,\\[3mm]
 |d_0|+|d_1|\leq 1.
 \end{array}
 \eeqtn
\item[\rm(iii)] For all $(d_0,d_1)\in D_T$ and $\frac{\varepsilon_0}{4}\leq |\theta|\leq \pi $, we have $u(\theta,d_0,d_1)=0$.
\end{enumerate}
\end{lemma}
{\textit{Proof}}: (i) and (ii) Since we have almost the same definition of the set $V_{K_0,A}$, and almost the same expression of initial data as in \cite{MZDuke97}, we refer the reader to Lemma 3.5 page 156 and Lemma 3.9 page 160 from \cite{MZDuke97}.\\ 
(iii) This follows by definition \eqref{initialq} of initial data for $s_0$ large enough.$\blacksquare$

\subsection{Details on the dynamics of equation (\ref{qqtilde})}
\label{reduction}
This subsection is dedicated to the proof of Lemmas \ref{transversality} and \ref{transcross}. They both follow from the understanding of the flow of equation (\ref{qqtilde}) in the set $S^*(t)$.\\ 
We proceed in two sections: We first prove Lemma \ref{transversality}, then Lemma \ref{transcross}.
\subsubsection{Reduction to a finite-dimensional problem}
Here we prove Lemma \ref{transversality}. Since the definition of $S^*(t)$ shows two different types of estimates, in the regions 
\pone$\;$and \ptwo, accordingly, we need two different approaches to handle those estimates:
\begin{itemize}
\item In \pone, we work in similarity variables (\ref{chauto}), in particular we crucially use the projection of equation (\ref{qqtilde}) with respect to the decomposition given in (\ref{decompq}).
\item In, \ptwo, we directly work in the variables $u(x,t)$, using standard parabolic estimates. 
\end{itemize}
\textbf{Part 1: Estimates in \pone.}

\medskip

\noindent In this part, we will show that 
\beqtn
\begin{array}{ll}
|q_2(s_{*})|\leq A^2 s^{-2}_{*} \log s_{*} -s_{*}^{3}, &              |q_-(y,s_{*})|\leq \frac{A}{2} s^{-2}_{*}(1+|y|^3),\\[3mm]
|q_e(y,s_{*})|\leq \frac{A^2}{2}s^{-2}_{*},&
\end{array}
\label{estiP1}
\eeqtn
where $s_*=-\log(T-t_*)$, $q$ is defined in (\ref{initialqqtilde}) and the notation is given in (\ref{decompq}), for a good choice of the parameters. In fact, this will follow from the projection of equation (\ref{qqtilde}) on the components $q_2$, $q_-$ and $q_e$ as we will see at the end of Part 1. Let us first give the behavior of those components in the following. 
\begin{prop}[Control of the null, negative and outer mode of equation (\ref{qqtilde})]\label{prop36}
For all $K_0>$, $\varepsilon_0>0$, there exists $\Acn\geq 1 $ such that for all $A\geq \Acn$ and $\etaa>0$, there exists $\scn(K_0,\varepsilon_0,A,\eta)$ such that the following holds for all $s_0\geq \scn$:\\
 Assume that for some $\tau\geq s_0$ and for all $s\in [\tau,\tau+\etaa]$,
 \[u(t)\in S^*(t)\mbox{, with }t=T-e^{-s}.\]
Then, the following holds for all $s\in [\tau,\tau+\etaa]$:

\begin{align*}
|q_{2}(s)|& \leq  \frac{\tau^2}{s^2}|q_2(\tau)|+\frac{CA(s-\tau)}{s^{3}},\\
\left\| \frac{q_-(s)}{1+|y|^3}\right\|_{L^{\infty}}&\leq C e^{-\frac{(s-\tau)}{2}} \left\| \frac{q_-(\tau)}{1+|y|^3}\right\|_{L^{\infty}}
+C\frac{e^{-(s-\tau)^2}\|q_e(\tau)\|_{L^\infty}}{s^{3/2}}+\frac{C(1+s-\tau)}{s^2},
\\
\|q_e(s)\|_{L^\infty}&\leq C e^{-\frac{(s-\tau)}{2}}\|q_e(\tau)\|_{L^\infty}+C \frac{e^{s-\tau}}{ s^{3/2}} \left\| \frac{q_-(\tau)}{1+|y|^3}\right\|_{L^\infty}+\frac{C(1+s-\tau)}{s^{1/2}}.
\end{align*}
\end{prop}
Let us first insist on the fact that the derivation of (\ref{estiP1}) follows from Proposition \ref{prop36}, exactly as in the real case treated in \cite{MZDuke97} (see pages 163 to 166 and 158 to 159 in \cite{MZDuke97} ). For that reason, we only focus in the following on the proof of Proposition \ref{prop36}.

\medskip

\subsubsection*{Proof of Proposition \ref{prop36}.}
\noindent The proof of Proposition \ref{prop36} consists in the projection of the equation (\ref{qqtilde}) on the different components of $q$ defined in (\ref{decompq}).\\
We note that the proof is already available from Lemma 3.13 page 167 and Lemma 3.8 page 158 from \cite{MZDuke97}, in the case of the standard heat equation in $\R^N$ without truncation terms.

\medskip

\noindent Since the equation satisfied by $q$ in (\ref{qqtilde}) shares the same linear part as the corresponding equation in \cite{MZDuke97}, the proof is similar to the argument in \cite{MZDuke97}, and the only novelty concerns the truncature term $F$ in (\ref{qqtilde}). For that reason, we only give the ideas here, focusing only on the new term $F$ and kindly ask the interested reader to look at Lemma 3.13 page 167 and Lemma 3.8 page 158 in \cite{MZDuke97} for the technical details.

\medskip

\noindent Let us first write equation \eqref{qqtilde} satisfied by $q$ in its Duhamel formulation,
\begin{eqnarray}
q(s)&=&K(s,\tau)q(\tau)+\int_{\tau}^{s}d\sigma K(s,\sigma)B(q(\sigma))+\int_{\tau}^{s}d\sigma K(s,\sigma)R(\sigma)\nonumber\\[3mm]
&+&\int_{\tau}^{s}d\sigma K(s,\sigma)(H+\pa_y G)(\sigma),
\end{eqnarray}
where $K$ is the fundamental solution of the operator $\Lg+V$.
We write $q=\alpha+\beta+\gamma+\delta+\tilde{\delta}$ with 
\begin{eqnarray}
\alpha(s)&=&K(s,\tau)q(\tau),\qquad\qquad\qquad\beta(s)=\int_{\tau}^{s}d\sigma K(s,\sigma)B(q(\sigma)),\nonumber\\[3mm]
\gamma(s)&=&\int_{\tau}^{s}d\sigma K(s,\sigma)R(\sigma),\qquad\quad\delta(s)=\int_{\tau}^{s}d\sigma K(s,\sigma)H(\sigma)\\[3mm] 
\tilde{\delta}(s)&=&\int_{\tau}^{s}d\sigma K(s,\sigma)\pa_yG(\sigma) \label{defdelta},\nonumber
\end{eqnarray}
where, for a function $F(y,\sigma)$, $K(s,\sigma)F(\sigma)$ is defined by
\[K(s,\sigma)F(\sigma)=\int_{\tau}^{s}d\sigma\int dx K(s,\sigma,y,x)F(x,\sigma).\]
We assume that $q(s)\in V_A(s)$ for each $s\in [\tau,\tau+\etaa]$. Clearly, proceeding as the derivation of Lemma 3.13 page 167 in \cite{MZDuke97}, Proposition \ref{prop36} follows from  the following:
\begin{lemma}[Projection of the Duhamel formulation] \label{lemduh}For all $K_0>0$, $\varepsilon_0>0$, there exists $\Adu\geq 1$ such that for all $A\geq \Adu$ and $\etaa>0$ there exists $\sdu(K_0,\varepsilon_0,A,\etaa)$, such that for all $s_0\geq \sdu$, if $0<\eta_0\leq 1$ and we assume that for some $\tau\geq s_0$ and for all $s\in[\tau,\tau+\etaa]$, $u(t)\in S^*(t)$ with $t=T-e^{-s}$, then
\begin{enumerate}
\item[\rm(i)] (Linear terms)
\beqtn
\begin{array}{lll}
|\alpha_2(s)|&\leq& \frac{\tau^2}{s^2}|q_2(\tau)|+\frac{C A  (s-\tau)}{s^3},\\[3mm]
\left\| \frac{\alpha_-(s)}{1+|y|^3}\right\|_{L^{\infty}}&\leq &C e^{-\frac{(s-\tau)}{2}} \left\| \frac{q_-(\tau)}{1+|y|^3}\right\|_{L^{\infty}}
+C\frac{e^{-(s-\tau)^2}\|q_e(\tau)\|_{L^\infty}}{s^{3/2}}+\frac{C}{s^2},\\[3mm]
\|\alpha_e(s)\|_{L^\infty}&\leq &Ce^{-\frac{(s-\tau)}{2}}\|q_e(\tau)\|_{L^\infty}+Ce^{s-\tau} s^{3/2} \left\| \frac{q_-(\tau)}{1+|y|^3}\right\|_{L^\infty}+\frac{C}{\sqrt{s}},
\end{array}
\eeqtn
\item[\rm(ii)]  (Nonlinear terms)
\begin{align*}
|\beta_2(s)|&\leq \frac{(s-\tau)}{s^3},&|\beta_-(y,s)|&\leq \frac{(s-\tau)}{s^2}(1+|y|^3),&\|\beta_e(s)\|_{L^\infty}&\leq \frac{(s-\tau)}{\sqrt s},\\[3mm]
|\delta_2(s)|&\leq C\frac{(s-\tau)}{s^3}, &|\delta_-(y,s)|&\leq C\frac{(s-\tau)}{s^2}(1+|y|^3),&\|\delta_e(s)\|_{L^\infty}&\leq C\frac{(s-\tau)}{\sqrt s},\\[3mm]
|\tilde{\delta}_2(s)|&\leq C\frac{(s-\tau)}{s^3},&|\tilde{\delta}_-(y,s)|&\leq C \frac{(s-\tau)}{s^2}(1+|y|^3),&\|\tilde{\delta}_e(s)\|_{L^\infty}&\leq C\frac{(s-\tau)}{\sqrt{s}}.
\end{align*}
\item[\rm(iii)]  (Source term)
\[|\gamma_2(s)|\leq C(s-\tau)s^{-3},\;|\gamma_-(y,s)|\leq C(s-\tau)(1+|y|^3)s^{-2},\;\|\gamma_e(s)\|_{L^\infty}\leq (s-\tau)s^{-1/2}.\]
\end{enumerate}
\end{lemma}
\textit{Proof:} We consider, $A\geq 1$, $\etaa>0$, 
and $s_0\ge \etaa$. 
The terms $\alpha$, $\beta$ and $\gamma$ are already present in the case of the real-valued semilinear heat equation, so we refer to Lemma 3.13 page 167 in \cite{MZDuke97} for the estimates involving them. Thus, we only focus on the new terms $\delta(y,s)$ and  $\tilde{\delta}(y,s)$.\\ Note that since $s_0\geq \etaa$, if we take $\tau\geq s_0$, then $\tau+\etaa\leq 2 \tau$ and if $\tau \leq \sigma\leq s\leq \tau +\etaa$, then
\beqtn
\frac{1}{2\tau}\leq \frac{1}{s}\leq \frac{1}{\sigma}\leq \frac{1}{\tau}.
\label{tausigma}
\eeqtn
Let us first derive the following bounds when $u(t)\in S^*(t)$:
\begin{lemma} 
 \label{VA}
For all $K_0>0$, $\varepsilon_0>0$, $A\geq 1$, there exists $s_0\geq \sdc (K_0,\varepsilon_0,A)$ such that if $s\geq \sdc$, $0<\eta_0\leq 1$ and we assume that $u(t)\in S^*(t)$ defined in Definition \ref{defVA}, where $t=T-e^{-s}$. Then, we have
 \[\begin{array}{l}
 (i)\; for \;all\; y \in \R,\;\; |q(y,s)|\leq C A^2 \frac{\log s}{s^2}(1+|y|^3),\\[3mm]
(ii)\;  \|q(s)\|_{L^\infty}\leq C\frac{A^2}{\sqrt{s}},\\[3mm]
(iii) \|W(s)\|_{L^{\infty}}\leq \kappa+2.

\end{array}
\]
 \end{lemma}
 \textit{Proof}. (i) and (ii): Since $u(t)\in S^*(t)$, it follows by definition that $q(s)\in V_{K_0,A}(s)$, where $s=-\log(T-t)$, therefore, the proof is the same as the corresponding part in \cite{MZDuke97}. See Proposition 3.7 page 157 in   \cite{MZDuke97} for details.\\
(iii)
 From (\ref{chauto}), (\ref{defw}) and (\ref{initialqqtilde}), we see that:
\begin{itemize}
\item If $|y|\leq \varepsilon_0 e^{s/2}$, then $W(y,s)=w(y,s)=\varphi(y,s)+q(y,s)$. Since $\|\varphi\|_{L^\infty}\leq \kappa +1$ from (\ref{deffi}), using (ii), we see that $\|W\|_{L^\infty}\leq\kappa+2$ for $s$ large enough that is for $T$ small enough.

\item If $|y|\geq \varepsilon_0 e^{s/2}$, then $W(y,s)=e^{-\frac{s}{p-1}}u(\theta e^{-s/2},t)$ with $|\theta|\geq \frac{\varepsilon_0}{2}$. By (ii) of Definition \ref{defVA}, we see that $|W(y,s)|\leq \eta_0 e^{-\frac{s}{p-1}}\leq \eta_0 T^{\frac{1}{p-1}}\leq 1$, if $\eta_0\leq 1$ and $T\leq 1$. 
\end{itemize}
This concludes the proof of Lemma \ref{VA}.$\blacksquare$

\medskip

\noindent Let us now recall from Bricmont and Kupiainen \cite{BKN94} the following estimates on $K(s,\sigma)$, the semigroup generated by $\Lg+V$.
\begin{lemma}[Properties of $K(s,\sigma)$]\label{lemK}
For all $s\geq \tau\geq1$, with $s\leq 2\tau$, we have the following:\\
{\rm (i)} for all $y$, $x\in\R$, we have,
\[|K(s,\sigma,y,x)|\leq C e^{(s-\sigma)\Lg}(y,x),\]
where $e^{\tt\Lg}$ is given explicitely by the Mehler's formula \cite{S79}
\beqtn
e^{\tt \Lg}(y,x)=\frac{e^\tt}{\sqrt{4\pi(1-e^{-\tt})}} \exp\left[-\frac{(ye^{-\tt/2}-x)^2}{4(1-e^{-\tt})}\right].
\label{defLG}
\eeqtn
{\rm (ii)} We have 
\beqtn
\displaystyle\left| \int K(s,\tau,y,x)(1+|x|^m)dx \right |\leq C\int e^{(s-\tau)\Lg}(y,x) (1+|x|^m)dx \leq e^{s-\tau} (1+|y|^m).
\eeqtn
{\rm (iii)} For all $g\in L^{\infty}$, such that $xg\in L^\infty$ 
\[
\|K(s,\tau)\pa_x g\|_{L^{\infty}}\leq C e^{s-\tau}\left \{\frac{\|g\|_{L^\infty}}{\sqrt{1-e^{-(s-\tau)}}}+\frac{(s-\tau)}{s}(1+s-\tau)\left ((1+e^{\frac{s-\tau}{2}})\|xg\|_{L^\infty}+ e^{\frac{s-\tau}{2}}\|g\|_{L^\infty}\right)\right\}.
\]
\end{lemma}
{\it Proof}.\\
(i) See page 181 in \cite{MZDuke97}\\
(ii) See Corollary 3.14 page 168 in \cite{MZDuke97}.\\
(iii) See Appendix \ref{appA}.
 $\blacksquare$
 
 \
 
 Now, with lemmas \ref{VA} and \ref{lemK} at hand, we are in position  to finish the proof of Lemma  \ref{lemduh}.  As we mentioned at the beginning of the proof, we only focus on the proof of the estimates on $\delta$ and $\tilde{\delta}$, and refer the readers to Lemma 3.13 page 167 in \cite{MZDuke97} for the estimate involving $\alpha$, $\beta$ and $\delta$.

 \medskip
 
 \noindent{\bf \small Estimates on $\delta$ defined in (\ref{defdelta})}:\\
 Consider $s\in[\tau,\tau+\etaa]$ and recall that $0<\eta_0\leq 1$. Since $u(t)\in S^*(t)$ with $t=T-e^{-s}$, we see from the definition  (\ref{chauto}) of $W$ that when
\beqtn
|y|\geq \frac{\varepsilon_0}{2}e^{s/2},\;	|W(y,s)|\leq \eta_0 e^{-\frac{s}{p-1}}\leq e^{-\frac{s}{p-1}}.
\label{ine1}
\eeqtn
Moreover by definition (\ref{Nchi}) of $\chi$, we see that 
\beqtn
\begin{array}{l}
|\pa_y\chi|\leq \frac{C}{\varepsilon_0}e^{-s/2}{\bf{I}}_{\varepsilon_0 e^{s/2}<|y|<2 \varepsilon_0 s^{s/2}},\\[3mm]
|y\pa_y\chi|\leq C {\bf{I}}_{\varepsilon_0 e^{s/2}<|y|<2 \varepsilon_0 s^{s/2}},\\[3mm]
\mbox{and }|\pa_s\chi|+(1+|y|)|\pa_y\chi|+|\pa_{y}^{2}\chi|+(\chi-\chi^p)\leq \frac{C}{\varepsilon_{0}^{2}} {\bf{I}}_{\varepsilon_0 e^{s/2}<|y|<2 \varepsilon_0 s^{s/2}}.
\end{array}
\label{ine2}
\eeqtn
Therefore, by definition (\ref{HG}), we see that 
\[\|H(s)\|_{L^\infty}\leq \frac{C}{\varepsilon_{0}^{2}} \eta_0 e^{-\frac{s}{p-1}}\leq \frac{C}{\varepsilon_{0}^{2}}  e^{-\frac{s}{p-1}}.\]
In particular, if $\tau\leq \sigma\leq s\leq \tau+\etaa$, we see from (\ref{tausigma}) that $\sigma \geq s/2$, hence 
\beqtn
\|H(\sigma)\|_{L^\infty}\leq \frac{C}{\varepsilon_{0}^{2}} e^{-\frac{\sigma}{p-1}}\leq \frac{C}{\varepsilon_{0}^{2}}e^{-\frac{s}{2(p-1)}}.
\label{estiH}
\eeqtn
Using Lemma \ref{lemK} and the  definition (\ref{defdelta}) of $\delta$, we write
\beqtn
\begin{array}{lll}
|\delta(y,s)|  &\leq&\displaystyle\int_{\tau}^{s} d\sigma\int_{\R}\left|K(s,\sigma,y,x)H(x,\sigma) \right|dx,\\[3mm]
&\leq&\displaystyle \int_{\tau}^{s} d\sigma\int_{\R} e^{(s-\sigma)\Lg}(y,x)\frac{C}{\varepsilon_{0}^{2}}e^{-\frac{s}{2(p-1)}} dx,\\[3mm]
&\leq &\displaystyle \frac{C}{\varepsilon_{0}^{2}}e^{-\frac{s}{2(p-1)}}\int_{\tau}^{s} d\sigma e^{(s-\sigma)},\\[3mm]
&\leq &\displaystyle \frac{C}{\varepsilon_{0}^{2}}e^{-\frac{s}{2(p-1)}} (s-\tau),\\[3mm]
&\leq &\displaystyle\frac{(s-\tau)}{s^2},
\label{esdelta}
\end{array}
\eeqtn
for $s$ large enough depending on $\eta_0$.\\
By definition of $q_m$, $q_-$ and $q_e$ for $m\leq 2$, we write
\beqtn
\begin{array}{lll}
|\delta_m(s)|&\leq&\left|\int_{\R}\chi(y,s)\delta(y,s)k_m(y)\rho(y)dy \right|\leq C\int_{\R}|\delta(y,s)|(1+|y|^2)\rho(y)dy\leq\frac{ C(s-\tau)}{s^{3}},\\[3mm]
|\delta_-(y,s)|&=&\left|\chi(y,s) \delta(y,s)-\sum_{i=0}^{2}\delta_i(s)k_i(y)\right|\leq(s-\tau)(1+|y|^3)\frac{C}{s^{2}}.\\[3mm]
|\delta_e(y,s)|&=&\left|(1-\chi(y,s))\delta (y,s)\right\|\leq(s-\tau)\frac{C}{\sqrt{s}},
\label{details}
\end{array}
\eeqtn
which are the desired estimations on $\delta$ in Lemma \ref{lemduh}.

\bigskip

 \noindent{\bf \small Estimates on $\tilde{\delta}$ defined in (\ref{defdelta})}:\\
Since for all $s\in [\tau,\tau+\etaa]$, $u(t)\in S^*(t)$, where $t=T-e^{-s}$, by assumption, using (\ref{ine1}) and (\ref{ine2}), we see that when $s_0\leq\tau\leq\sigma\leq\tau\leq \tau+\etaa$, we have $\sigma\geq s/2$, hence
\beqtn\label{estiG1}
\|G(\sigma)\|_{L^\infty}\leq \frac{C\eta_0}{\varepsilon_0} e^{-\frac{(p+1)\sigma}{2(p-1)}}\leq \frac{C}{\varepsilon_0} e^{-\frac{(p+1)s}{4(p-1)}},
\eeqtn
\beqtn\label{estiG}
\|x G(\sigma)\|_{L^\infty}\leq C\eta_0 e^{-\frac{\sigma}{p-1}}\leq C e^{-\frac{s}{2(p-1)}},
\eeqtn
where $G$ is defined by (\ref{HG}), remember that $\eta_0\leq1$.\\
Using (iii) of Lemma \ref{lemK}, with $g=G(\sigma)$, we obtain
 \begin{eqnarray*}
&&\|K(s,\sigma)\pa_x G\|_{L^{\infty}}\leq C e^{s-\sigma}\bigg \{\frac{e^{-\frac{(p+1)s}{4(p-1)}}}{\sqrt{1-e^{-(s-\sigma)}}}\\[3mm]
&&\qquad\qquad\qquad\qquad\qquad+\frac{(s-\sigma)}{s}(1+s-\sigma)\bigg(e^{-\frac{s}{2(p-1)}}(1+e^{\frac{s-\sigma}{2}})+e^{\frac{s-\sigma}{2}}e^{-\frac{(p+1)s}{4(p-1)}}\bigg)\bigg\}.
\end{eqnarray*}
Integrating in time, we get rid of the square rest term in the denominator, and see that 
\[\left|\int_{\tau}^{s}K(s,\sigma)\pa_x G(x,\sigma)d\sigma\right|\leq C e^{s-\tau}(s-\tau)e^{-\frac{s}{2(p-1)}}\left(e^{-\frac s 4}+ \frac{(s-\tau)}{s}(1+s-\tau)(1+e^{(s-\tau)/2})\right).\]
Proceeding as for $\delta$ and using the fact that $0\leq s-\tau \leq \etaa  $, we get 
\[\sum_{m=0}^{2}|\tilde{\delta}_m(s)|+\left\|\frac{\tilde{\delta}_m(s)}{1+|y|^3}\right\|_{L^\infty}+\|\tilde{\delta}_e(s)\|_{L^\infty}\leq \frac{C(s-\tau)}{s^3},\]
for $s_0$ large enough, which gives the desired estimates on $\tilde{\delta}$ in Lemma \ref{lemduh}. Since the estimate notes for $\alpha$, $\beta$ and $\gamma$
defined in (\ref{defdelta}) follows exactly as in Lemma 3.13 p 167 in \cite{MZDuke97}, this concludes the proof of Lemma \ref{lemduh}.$\blacksquare$\\

\noindent {\bf \small Conclusion of Part 1 and choice of parameters:} Proceeding exactly as in \cite{MZDuke97}, page 157, we derive estimate \eqref{estiP1} from Proposition \ref{prop36}. This is possible for any $K_0>0$, $\varepsilon_0>0$, $A\geq \Acn$, for some $\Acn(K_0,\varepsilon_0)\geq 1$ and $s_0\geq \scn$ for some $\scn(K_0,\varepsilon_0,A)$.

\medskip

\noindent Since Proposition \ref{prop36} follows directly from Lemma  \ref{lemduh}, this ends the proof of Proposition \ref{prop36}, as we mentioned right before the statement of Lemma \ref{lemduh} in the same way as in Proposition 3.11 page 161 in \cite{MZDuke97}.$\blacksquare$
\label{pagepart1}

\medskip

\noindent \textbf{Part 2: Estimates in \ptwo.}

\medskip

\noindent The aim of this part is to show that
\beqtn
\mbox{if  }\frac{\varepsilon_0}{2}\leq |\theta|\leq \pi\mbox{, then }|u(\theta,t_*)|\leq \frac{\eta_0}{2},
\label{estitetaeps}
\eeqtn
 provided the parameters satisfy some conditions. 
We proceed in 3 steps:
\begin{itemize}
\item In Step 1, we derive better bounds on the solution $u(\theta,t)$ in the intermediate region
\beqtn
K_0\sqrt{(T-t)\left |\log(T-t)\right |}\leq |\theta|\leq \frac{\varepsilon_0}{2}.
\label{esti1}
\eeqtn
\item In Step 2, we introduce a parabolic estimate on the solution in the region \ptwo.
\item Finally, in Step 3, we combine the previous steps to show (\ref{estitetaeps}).
 \end{itemize}
 \textbf{Step 1: Improved estimates in the intermediate region.}\\
Here, we refine the estimates on the solution in the region (\ref{esti1}).
In fact, we have from item (iii) of Lemma \ref{VA}
\beqtn
\forall t\in[0,t_*],\;\forall \theta \in \R,\;  |u(\,t)|\leq C (T-t)^{-\frac{1}{p-1}},
\label{estitaux}
\eeqtn
valid in particular in the region (\ref{esti1}). This bound is not satisfactory, since it goes to infinity as $t\to T$. In order to refine it, given a small $\theta$, we use this bound when $t=t_0(\theta)$ defined by
\beqtn
\label{estiteta0}
|\theta|=K_0\sqrt{(T-t_0(\theta))|\log(T-t_0(\theta))|},
\eeqtn
to see that the solution is in fact flat at that time. Then, advancing the PDE (\ref{eqtheta}), we see that the solution remains flat for later times. More precisely, we claim the following:
\begin{lemma}[Flatness of the solution in the intermediate region in (\ref{esti1})]
There exists $\zeta_0>0$ such that for all $K_0>0$, $\varepsilon_0>0$, $A\geq 1$, there exists $\sfl(K_0,\varepsilon_0,A)$, such that if $s_0\geq \sfl$ and $0<\eta_0\leq 1$, then, 
\[\forall t_0(\theta)\leq t \leq t_*,\;\; \left|\frac{u(\theta,t)}{u^*(\theta)}-\frac{U_{K_0}(\theta)}{U_{K_0}(1)} \right|  \leq  \frac{C}{|\log\theta|^{\zeta_0}},\]
where $u^*$ is defined in \eqref{defue} and
\beqtn
U_{K_0}(\tau)=\kappa \left((1-\tau)+\frac{(p-1)K_{0}^{2}}{4p}\right)^{-1/(p-1)}
\label{UK0}
\eeqtn
In particular, $|u(\theta,t)|\leq 2 |u^*(\theta)|$. 

\label{flatnes}
\end{lemma}
\textit{Proof: }We argue as in Masmoudi and Zaag \cite{MZ07}. If $\theta_0\neq 0$ is small enough, we introduce for all $(\xi,\tau)\in \R\times [-\frac{t_0(\theta_0)}{T-t_0(\theta_0)},\tau_*)$, with $\tau_*=\frac{t_*-t_0(\theta_0)}{T-t_0(\theta_0)}$
\begin{equation}
\label{defV}U(\theta,\xi,\tau)=(T-t_0(\theta_0))^{1/(p-1)}u(\theta,t),
\end{equation}
where
\begin{equation}\theta=\theta_0+\xi\sqrt{T-t_0(\theta_0)},\; t=t_0(\theta_0)+\tau(T-t_0(\theta_0)),
\end{equation}
and $t_0(\theta_0)$ is uniquely defined by
\beqtn
\label{estiteta00}
|\theta_0|=K_0\sqrt{(T-t_0(\theta_0))|\log(T-t_0(\theta_0))|},
\eeqtn
From the invariance of problem (\ref{eqtheta}) under dilation, $U(\theta_0,\xi,\tau)$ is also a solution of (\ref{eqtheta}) on its domain. Since $u(t)\in S^*(t)$, from \eqref{us}, using the definition of $S^*(t)$, Lemma \ref{VA} together with (\ref{estiteta0}) and (\ref{defV}), we have

\[\begin{array}{ll}
\displaystyle\sup_{|\xi|<2|\log(T-t_0(\theta_0))|^{1/4}}\left |U(\theta_0,\xi,0)-f(K_0) \right|&\leq\displaystyle \frac{K_0}{2p\mathcal{S}_0}+\|q(\mathcal{S}_0)\|_{L^\infty}\\[3mm]
&\leq\displaystyle\frac{C}{\mathcal{S}_0}+\frac{C A^2}{\mathcal{S}_{0}^{1/2}}\\[3mm]
&\leq\displaystyle\frac{C}{\mathcal{S}_{0}^{1/4}},
\end{array}
\]
with $\mathcal{S}_0=\mathcal{S}_0(\theta_0)=-\log(T-t_0(\theta_0))$ and $f$ is defined by \eqref{deff}. Provided that $s_0(=-\log T)$ is large enough.\\
Using the continuity with respect to initial data for problem (\ref{eqtheta}), associated to a space-localization in the ball $B(0,|\log(T-t_0(\theta_0))|^{1/4})$, we show as in Section 4 of \cite{ZAIHPANL98} that
\beqtn
\begin{array}{l}
\sup_{|\xi|\leq |\log(T-t_0(\theta_0))|^{1/4},\;0\leq\tau<\tau_*}|U(\theta_0,\xi,\tau)-U_{K_0}(\tau)|\leq \frac{C}{\mathcal{S}_{0}^{\zeta_0}},
\end{array}
\label{estiteta01}
\eeqtn
$U_{K_0}(\tau)$ given by \eqref{UK0} is the solution of the PDE (\ref{eqtheta}) with constant initial data $f(K_0)$. Since $U_{K_0}(\tau)\leq U_{K_0}(1)=\kappa\left(\frac{(p-1)K_{0}^{2}}{4p}\right)^{-1/(p-1)}$ and we have from \eqref{estiteta00}
\beqtn
\log(T-t_0(\theta_0))\sim 2\log\theta_0\mbox{ and }
\left(T-t_0(\theta_0)\right)\sim \frac{\theta_{0}^{2}}{2K_{0}^{2}|\log(\theta_0)|}\mbox{ as }\theta_0\to 0,
\label{estit0}
\eeqtn
this yields $(T-t_0(\theta_0))^{1/(p-1)}\sim\displaystyle \frac{U_{K_0}(1)}{u^*(\theta_0)}$, by definition \eqref{defue} of $u^*$.
We obtain the desired conclusion from (\ref{defV}) and (\ref{estiteta01}). This ends the proof of Lemma \ref{flatnes}.$\blacksquare$

\medskip

\noindent\textbf{Step 2: A parabolic estimate in Region \ptwo}\\
We recall from the definition \ref{defVA} of $S^*(t)$ that
\[\forall \theta\in \R\mbox{ such that }\frac{\varepsilon_0}{2}\leq |\theta|\leq \pi,\;\;|u(\theta,t)|\leq \eta_0.\]
Here, we will obtain a parabolic estimate on the solution in \ptwo. More precisely, we claim the following: 
\begin{prop}[A parabolic estimate in \ptwo]\label{parabolicesti}
For all $\varepsilon>0$, $\varepsilon_0>0$, $\sigma_1\geq 0$, $\exists \,T_4(\varepsilon,\varepsilon_0,\sigma_1)\geq 0$, such that for all $\bar{t}\leq T_4$, if $u$ a  periodic solution of 
\[\pa_t u=\pa_{\theta}^{2}u+|u|^{p-1}u\mbox{ for all } \theta\in \mathbb{S},\;\; t\in [0,\bar{t}],\]
which satisfies:
\begin{enumerate}
\item[\rm(i)] for $|\theta|\in [\frac{\varepsilon_0}{4},\frac{\varepsilon_0}{2}],\;\;|u(\theta,t)|\leq \sigma_1$.
\item[\rm(ii)] for $\frac{\varepsilon_0}{4}\leq |\theta|\leq \pi$, $u(\theta,0)=0$.
\end{enumerate}
Then, for all $t\in[0,\bar t]$, for all $\frac{\varepsilon_0}{2}\leq|\theta|\leq\pi$,
\[|u(\theta,t)|\leq \varepsilon.\]
\end{prop}

\textit{Proof}:
Consider $\bar{u}$ defined in (\ref{defbaru}), which satisfies equation (\ref{eqbaru}), recalled here, after a trivial chain rule to transform the $\pa_\theta u$ term: 
\[\forall t\in[0,\bar{t}],\;\forall \theta\in \R,\;\pa_t \bar{u}=\pa_{\theta}^{2}\bar{u}+|u|^{p-1}\bar{u}-2\pa_\theta(\bar{\chi}'u)+\bar{\chi}''u.\]
Therefore, since $\bar{u}(\theta,0) \equiv 0$, we write
\[
\|\bar{u}(t)\|_{L^\infty}\leq \displaystyle\int_{0}^{t} \left| S(t-t')\left[|u|^{p-1}\textbf{I}_{|\theta|\geq \frac{\varepsilon_0}{4}}\bar{u}-2\pa_\theta\left(\bar{\chi} ' u\textbf{I}_{|\theta|\geq \frac{\varepsilon_0}{4}}\right)+\bar{\chi}^{''}u(t')\textbf{I}_{|\theta|\geq \frac{\varepsilon_0}{4}}\right]\right|dt' , 
\]

where $S(t)$ is the heat kernel.\\
Since $\bar{\chi}^{'}$ and $\bar{\chi}^{''}$ are supported by $\{\frac{\varepsilon_0}{4}\leq |\theta|\leq \frac{\varepsilon_0}{2}\}$ and satisfy $|\bar{\chi}^{'}|\leq C/\varepsilon_0$,  $|\bar{\chi}^{''}|\leq C/\varepsilon_{0}^{2}$ and using parabolic regularity, we write
\[
\begin{array}{ll}
\|\bar{u}(t)\|_{L^\infty}&\leq \sigma_{1}^{p-1}\int_{0}^{t}\|\bar{u}(t')\|dt'+\frac{C\sigma_1}{\varepsilon_0}\int_{0}^{t}\frac{dt'}{\sqrt{t-t'}}+\frac{C\sigma_1}{\varepsilon_{0}^{2}}\int_{0}^{t}dt'\\[3mm]
&\leq \sigma_{1}^{p-1}\int_{0}^{t}\|\bar{u}(t')\|dt'+\frac{C\sigma_1}{\varepsilon_0}\sqrt{\bar{t}}+\frac{C\sigma_1}{\varepsilon_{0}^{2}}\bar{t}.
\end{array}
\]
If $\bar{t}<1$, by Gronwall estimate, this implies that 
\[ \|\bar{u}(t)\|_{L^\infty}\leq Ce^{\sigma_{1}^{p-1}}\left(\frac{\sigma_1}{\varepsilon_0}\sqrt{\bar{t}}+\frac{\sigma_1}{\varepsilon_{0}^{2}}\bar{t}\right).\]
Taking $\bar{t}$ small enough, we can obtain
\[\forall t \in [0,\bar{t}],\;\;\|\bar{u}(t)\|_{L^\infty}\leq \varepsilon.\]
Since $u=\bar{u}$ for all $\frac{\varepsilon_0}{2}\leq |\theta|\leq \pi$ by definition (\ref{defbaru}), this concludes the proof of Proposition \ref{parabolicesti}. $\blacksquare$\\

\noindent\textbf{Step 3: Proof of the improvement in (\ref{estitetaeps})}\\
Here, we use Step 1 and Step 2 to prove (\ref{estitetaeps}), for a suitable choice of parameters.\\
Let us consider $K_0>0$, and $\delta_0(K_0)>0$ defined in Lemma \ref{flatnes}. Then, we consider $\varepsilon_0\leq 2\delta_0$, $A\geq 1$, $0<\eta_0\leq 1$, and
\[s_0\geq \spi \equiv \max\{\sfl(K_0,\varepsilon_0,A),\sdi(K_0,\varepsilon_0,A),-\log\left(T_4\left(\frac{\eta_0}{2},\varepsilon_0,2 |u^*(\frac{\eta_0}{4})|\right)\right) \},\]
where the different constants are defined in Lemmas \ref{prop37} and \ref{flatnes} and Proposition \ref{parabolicesti}.\\
Applying Lemma \ref{flatnes}, we see that
\[\forall \,|\theta|\leq \delta_0,\; A\geq 1,\;\forall t\in[0,t_*],\; |u(\theta,t)|\leq 2 |u^*(\theta)|.\]
In particular,

\[\forall \,\frac{\varepsilon_0}{4}\leq |\theta|\leq \frac{\varepsilon_0}{2}\leq \delta_0,\; \forall t\in[0,t_*],\; |u(\theta,t)|\leq 2 |u^*(\frac{\varepsilon_0}{4})|.\]
Using item (iii) of Lemma \ref{prop37}, we see that $\forall \frac{\varepsilon_0}{4}\leq |\theta|\leq \pi$, $u(\theta,0)=0$.\\
Therefore Proposition \ref{parabolicesti} applies with $\varepsilon= \frac{\eta_0}{2}$ and $\sigma_1=2u^*( \frac{\varepsilon_0}{4})$ and we see that
\beqtn
\forall \frac{\varepsilon_0}{2}\leq |\theta|\leq \pi,\; \forall t\in[0,t_*],\; |u(\theta,t)|\leq \frac{\eta_0}{2}.
\label{inetheeta}
\eeqtn
and estimate (\ref{estitetaeps}) holds.

\medskip

\noindent{\bf \small Conclusion of Part 2 and choice of parameters:} From \eqref{inetheeta}, we see that (\ref{estitetaeps}) holds for any $K_0>0$, $\varepsilon_0\leq 2 \delta_0(K_0)$, $A\geq 1$, $0<\eta_0\leq 1$ and $s_0\geq \spi(K_0,\varepsilon_0,A,\eta_0)$.
\label{pagepart2}

\medskip

{\it Conclusion of the proof of Lemma \ref{transversality}:} From the conclusion of Parts 1 and 2, if we take $K_0>0$, $\varepsilon_0\leq 2\delta_0(K_0)$, $A\geq \Acp (K_0,\varepsilon_0)$, $0<\eta_0\leq 1$ and 
\[s_0\geq \max\{\scn(K_0,\varepsilon_0,A),\spi (K_0,\varepsilon_0,A,\eta_0)\},\]
then we see that \eqref{estiP1} and (\ref{estitetaeps}) holds at $t=t_*$.\\
 Recalling that $u(t_*)\in \pa S(t_*)$ by \eqref{336}, we see from Definition \ref{defVA} of $S^*$ that only one of the components $q_0(s_*)$ or $q_1(s_*)$ may touch the boundary of $[-\frac{A}{s_{*}^{2}},\frac{A}{s_{*}^{2}}]$. This concludes the proof of Lemma \ref{transversality}.$\blacksquare$
\label{subsect331}
\subsubsection{Transverse crossing on $V_{K_0,A}(s)$}
We prove Lemma \ref{transcross} here. The key estimate is to prove the following differential inequality on $q_m$ for $m=0,1$:
\beqtn
\forall s\in [s_0,s_*],\;\; \left|q_{m}^{'}(s)-(1-\frac m2)q_m(s)\right|\leq \frac{C}{s^2},
\label{inetrans}
\eeqtn
provided $s_0\geq \str(K_0,\varepsilon_0,A)$ and $0<\eta_0\leq 1$, for some large enough $\str$. Indeed if (\ref{qmqtildem}) holds, say $q_m(s_*)=\frac{\omega A}{s_{*}^{2}}$ for $m=0,1$ and $\omega=\pm 1$, then, we see that 
\[\omega q_{m}^{'}(s)\geq  (1-\frac m2)\frac{A}{s_{*}^{2}}-\frac{C}{s^2}\geq (1-\frac m2)\frac{A}{2s_{*}^{2}},\]
assuming that $A$ is large enough, which yields the conclusion of Lemma \ref{transcross}, assuming that (\ref{inetrans}) holds.

\medskip

\noindent Let us briefly justify (\ref{inetrans}). Multiplying equation (\ref{qqtilde}) by $h_m(y)\chi_1(y,s)$, defined in (\ref{chi1}), we obtain the following estimate
\[(\pa_s q)_m=\displaystyle ({\cal{L}}q)_m+(Vq)_m+B_m+R_m+F_m.\]
From straightforward estimates, already used for the standard heat equation in $\R^N$ considered in \cite{MZDuke97} (see Lemma 3.8 page 158 there), we know that
\[\left|(\pa_s q)_m-q^{'}_{m}\right|+\left|({\cal{L}}q)_{m}-(1-\frac m2)q_m\right|+\left|(Vq)_m+B_m+R_m\right|\leq \frac{C}{s^2}.\]

It remains only to treat the new term $F_m$. In fact from (\ref{estiH}), we see that 
\[\|H(s)\|_{L^\infty}\leq C e^{-\frac{s}{2(p-1)}}\leq \frac{C}{s^2}, \mbox{ for $s_0$ large enough, provided that $0<\eta_0\leq 1$}. \]
Then integrating by parts and using \eqref{ine1}, \eqref{ine2}, \eqref{estiG1} and \eqref{estiG}, we get
\[\left |\left(\pa_yG \right)_m\right |\leq \frac{C}{s^2}\mbox{ for $s_0$ large enough,}\]
and we obtain the following
\[|F_m(s)|\leq \frac{C}{s^2}.\]
This concludes the proof of (\ref{inetrans}) and Lemma \ref{transcross} too. $\blacksquare$

\section{Proof of Theorem \ref{theorem1}}
We prove Theorem \ref{theorem1} in this section. We will first derive (ii) from Section 3, then we will prove (i) and (iii). \\
Let us fix $K_0>0$, $\varepsilon_0>0$, $A>0$, $0<\eta_0\leq 1$ and $T>0$ so that Proposition \ref{Newprop} as well as all the statements of Section 3 apply hence, for some $d_0,d_1\in\R^2$, equation \eqref{eqtheta} with initial data given by \eqref{initialq} has a solution $u(\theta,t)$ such that
\beqtn
T=t_*(d_0,d_1),\;\forall t\in [0,T), u(t)\in S^*(K_0,\varepsilon_0,A,\eta_0,T,t).
\label{Ttstar}
\eeqtn
(Note the fact that $t_*(d_0,d_1)=T$ follows from the conclusion of the topology argument given by \eqref{conpro}.)\\
Applying item (ii) of Lemma  \ref{VA}, we see that
\[\forall y \in \R,\;\forall s\geq -\log T,\; |q(y,s)|\leq \frac{C A^2}{\sqrt{s}}.\]
By defintions \eqref{defw}, \eqref{initialqqtilde} and \eqref{deffi}, we see that
\[\forall s\geq -\log T,\;\forall |y|\leq \varepsilon_0 e^{s/2},\;\; \left| W(y,s)-f\left(\frac{y}{\sqrt{s}}\right)\right|\leq \frac{CA^2}{\sqrt{s}}+\frac C s\]
By definition \eqref{chauto} of $W$, we see that 
\[\forall t\in[0,T),\;\forall |\theta|\leq \varepsilon_0,\;\; \left|(T-t)^{1/(p-1)}u(\theta,t)-f\left(\frac{\theta}{\sqrt{(T-t)|\log(T-t)|}} \right)\right|\leq \frac{C(A)}{\sqrt{|\log(T-t)|}}.\] 
Since $u$ is $2\pi-$periodic estimate (\ref{profilev}) holds.

\medskip

\noindent(i) If $\theta_0=2k\pi,\; k\in \Z$, then we see from (\ref{profilev}) that $|u(0,t)| \sim\kappa (T-t)^{-1/(p-1)}$ as $t\to T$, where $\kappa$ is defined in \eqref{defkappa}. Hence $u$ blows up at time $T$ at $\theta_0=2k\pi,\; k\in\Z$.

\medskip

\noindent It remains to prove that any $\theta_0\neq 2k\pi$ is not a blow-up point.\\ 
From periodicity, we may assume that $-\pi \leq\theta_0\leq \pi$.\\ 
Since, we know from item (ii) in the Definition \ref{defVA}, that if $\frac{\varepsilon_0}{2}\leq |\theta|\leq \pi$, and $0\leq t\leq T$, $|u(\theta,t)|\leq \eta_0$, it follows that $\theta_0$ is not a blow-up point, provided 
\[\frac{\varepsilon_0}{2}\leq |\theta_0|\leq \pi.\]
Now, if $0<|\theta_0|\leq \frac{\varepsilon_0}{2}$, the following result from Giga and Kohn \cite{GKCPAM89} allows us to conclude.

\begin{prop}[Giga and Kohn - No blow-up under the ODE threshold] For all $C_0>0$, there is $\eta_0>0$ such that if $v(\xi,\tau)$ solves
\[\left |  v_t -\Delta v\right |\leq C_0 (1+|v|^p)\]
and satisfies
\[|v(\xi,\tau)|\leq \eta_0(T-t)^{-1/(p-1)}\]
for all $(\xi,\tau)\in B(a,r)\times[T-r^2,T)$ for some $a\in \R$ and $r>0$, then $v$ does not blow up at ($a,T$).
\label{propositionGK}
\end{prop}

\textit{Proof: } See Theorem 2.1 page 850 in \cite{GKCPAM89}. $\blacksquare$\\
Indeed, since $|\theta_0|\leq \frac{\varepsilon_0}{2}$, it follows from (\ref{profilev}) that
\[\sup_{|\theta-\theta_0|\leq |\theta_0|/2}(T-t)^{\frac{1}{p-1}}|u(\theta,t)|\leq \left | f\left(\frac{|\theta_0|/2}{\sqrt{(T-t)|\log(T-t)|}}\right)\right |+\frac{C}{\sqrt{|\log(T-t)|}}\to 0\]
as $t\to T$. Therefore, applying Proposition \ref{propositionGK}, we see that $\theta_0$ is not a blow-up point of $u$. This concludes the proof of (i) of Theorem \ref{theorem1}.

\medskip

(iii) Arguing as Merle did in \cite{FMCPAM92}, we derive the existence of a blow-up profile $u(\theta,T)\in C^2(\R\setminus \{2k\pi,k\in\Z\})$ such that $u(\theta,t)\to u(\theta,T)$ as $t\to T$, uniformly on compact sets of $\R\setminus \{2k\pi,k\in\Z\}$. The profile $u(\theta,t)$ is not defined at the origin. In the following, we would like to find its equivalent as $\theta\to 2k\pi$, for any $k\in\Z$ and show that it is in fact singular at $\theta_0=2k\pi$.\\ 
From periodicity it is enough to take $\theta_0=0$. Since $t_*(d_0,d_1)=T$ from \eqref{Ttstar}, applying Lemma \ref{flatnes} and making $t\to T$. We see that
\[\left|\displaystyle \frac{u(\theta,T)}{u^*(\theta)}-1\right|\leq \frac{C}{|\log\theta|^{\zeta_0}}.\]
Making $\theta\to 0$, we get the desired estimate in item (iii). This concludes the proof of Theorem \ref{theorem1}. $\blacksquare$


\appendix
 
\section{Proof of (iii) of Lemma \ref{lemK}}\label{appA}
This part relies mainly on the understanding of the behavior of the kernel $K(s,\sigma,y,x)$. This behavior follows from a perturbation method around $e^{(s-\sigma)\Lg}(y,x)$.\\
Since $\Lg$ is conjugated to the harmonic oscillator $e^{-x^2/8}\Lg e^{x^2/8}=\pa^2-(x^2/16)+(1/4)+1$, We use the definition of $K$ as the semigroup generated by $\Lg +V$ defined in \eqref{OperatorL} and give a Feynman-Kac representation for $K$:
\beqtn
K(s,\sigma,y,x)=e^{(s-\sigma)\Lg}(y,x)E(y,x),
\label{defK}
\eeqtn
where 
\beqtn
E(y,x)=\int d\mu_{yx}^{s-\sigma}(\omega)e^{\int_{0}^{(s-\sigma)}V(\omega(\tau),\sigma+\tau)},
\label{defE}
\eeqtn
and $d\mu_{yx}^{s-\sigma}$ is the oscillator measure on the continuous paths $\omega:[0,s-\sigma]\to \R$ with $\omega(0)=x$, $\omega(s-\sigma)=y$, i.e, the Gaussian probability measure with covariance kernel 

\[\Gamma(\tau,\tau')=\omega_0(\tau)\omega_0(\tau')+2\left(e^{-\frac 12 |\tau-\tau'|}-e^{-\frac 12 |\tau+\tau'|}+e^{-\frac 12 |2(s-\sigma)+\tau-\tau'|}-e^{-\frac 12 |2(s-\sigma)-\tau-\tau'|}\right),\]
which yields $\displaystyle \int d\mu_{yx}^{s-\sigma}\omega(\tau)=\omega_0(\tau)$ with
\[\omega_0(\tau)=\left(\sinh(\frac{s-\sigma}{2})\right)^{-1}\left(y \sinh \frac\tau 2+x\sinh \frac{s-\sigma-\tau}{2}\right).\]
Consider $1\leq \tau\leq s\leq 2\tau$. From Bricmont and Kupiainen \cite{BKN94}, Lemma 6 page 555 and Merle and Zaag \cite{MZDuke97} pages 183-184, we have the following estimates:
\begin{cl} If $1\leq \tau\leq s$ with $s\leq 2 \tau$, then
\[0\leq E(y,x)\leq C,\]
\[|\pa_xE(y,x)|\leq \frac{C}{s}(s-\sigma)(1+s-\sigma)(|y|+|x|).\]
\label{cl1}
\end{cl}

\medskip

Consider $g\in L^\infty$ such that $xg \in L^\infty$. By (\ref{defK}), (\ref{defE}) and integration by parts, we obtain
\beqtn
\begin{array}{ll}
\displaystyle K(s,\sigma)(\pa_x g)(y)&=\int_{\R}dxK(s,\sigma,y,x)\pa_x g(x,\sigma)=\displaystyle\int_{\R}dx e^{(s-\sigma)\Lg}(y,x)E(y,x) \pa_x g(x)\\[3mm]
&=\displaystyle-\int_{\R}dx\pa_x e^{(s-\sigma)\Lg}(y,x)E(y,x) g(x)-\int_{\R}dx e^{(s-\sigma)\Lg} \pa_xE(y,x)g(x)\\[3mm]
&=I_1+I_2.
\end{array}
\eeqtn

By (\ref{defLG}), we have
\begin{eqnarray*}
|\pa_x e^{(s-\sigma)\Lg}(y,x)|&=&\left|\frac{e^{s-\sigma}}{\sqrt{4\pi (1-e^{-(s-\sigma)})}}\frac{-(x-ye^{-\frac{(s-\sigma)}{2}})}{2 (1-e^{-(s-\sigma)})} \exp\left[-\frac{(x-ye^{-(s-\sigma)/2})^2}{4(1-e^{-(s-\sigma)})}\right]\right|\\[3mm]
&=&\displaystyle C\frac{|z|e^{-z^2}}{\sqrt{1-e^{-(s-\sigma)}}} \frac{e^{s-\sigma}}{\sqrt{1-e^{-(s-\sigma)}}}
\end{eqnarray*}
where  $z=\frac{x-ye^{-(s-\sigma)/2}}{\sqrt{4 (1-e^{-(s-\sigma)})}}$. 
Using Claim \ref{cl1}, we see that
\begin{eqnarray}
|I_1|\leq\displaystyle  C\|g\|_{L^\infty}\int_{\R}dx \frac{|z|e^{-z^2}}{\sqrt{1-e^{-(s-\sigma)}}} \frac{e^{(s-\sigma)}}{\sqrt{1-e^{-(s-\sigma)}}}&=&C\|g\|_{L^\infty}\int_{\R}dz  |z|e^{-z^2}\frac{e^{(s-\sigma)}}{\sqrt{1-e^{-(s-\sigma)}}} \nonumber\\[3mm]
&\leq& C \displaystyle\|g\|_{L^\infty}\frac{e^{(s-\sigma)}}{\sqrt{1-e^{-(s-\sigma)}}},
\end{eqnarray}
and
\begin{eqnarray}
|I_2|&\leq&  \int_{\R} dx e^{(s-\sigma)\Lg}(y,x)|g(x)|\frac{C(s-\sigma)}{s}(1+s-\sigma)(|y|+|x|),\nonumber\\[3mm]
&\leq&  \frac{C(s-\sigma)}{s}(1+s-\sigma)\int_{\R} dx e^{(s-\sigma)\Lg}(y,x)|g(x)|\left(|z|e^{(s-\sigma)/2}\sqrt{1-e^{-(s-\sigma)}}+|x|(1+e^{(s-\sigma)/2})\right),\nonumber\\[3mm]
&=&\frac{C(s-\sigma)}{s}(1+s-\sigma)(J_1+J_2),
\end{eqnarray}
with  $z=\frac{x-ye^{-(s-\sigma)/2}}{\sqrt{4 (1-e^{-(s-\sigma)})}}$.
Moreover
\[
\begin{array}{lll}
J_1&\leq&\displaystyle \|g\|_{L^\infty}e^{(s-\sigma)/2}\sqrt{1-e^{-(s-\sigma)}}\int_{\R} dx e^{(s-\sigma)\Lg}(y,x)|z|,\\[3mm]
&=&\displaystyle \|g\|_{L^\infty}e^{3 (s-\sigma)/2}\left(1-e^{-(s-\sigma)}\right) \int dz e^{-z^2}|z|,\\[3mm]
&\leq&\displaystyle C e^{3(s-\sigma)/2}\|g\|_{L^\infty}.
\end{array}
\]
Furthermore using item (i) of Lemma \ref{lemK}, we write
\[
\begin{array}{lll}
J_2&=&\displaystyle 2\int_{\R} dx e^{(s-\sigma)\Lg}(y,x)|x g(x)|(1+e^{(s-\sigma)/2}),\\[3mm]
&\leq&\displaystyle  C \|xg\|_{L^\infty}e^{s-\sigma}(1+e^{(s-\sigma)/2}).
\end{array}
\]
Which gives the desired estimates thanks to item (i) of Lemma \ref{lemK}.$\blacksquare$


\end{document}